\documentclass[11pt,intlimits,reqno]{amsart}

\usepackage[hypertex]{hyperref}
\usepackage{graphicx}
\usepackage{amsfonts,amssymb}
\usepackage{amsthm}
\usepackage{txfonts}

\usepackage{a4wide}

\newtheorem{neu}{}[section]
\newtheorem{Cor}[neu]{Corollary}
\newtheorem*{Cor*}{Corollary}
\newtheorem{Thm}[neu]{Theorem}
\newtheorem*{Thm*}{Theorem}
\newtheorem{Prop}[neu]{Proposition}
\newtheorem*{Prop*}{Proposition}
\theoremstyle{definition}
\newtheorem{Lemma}[neu]{Lemma}
\newtheorem*{Rmk*}{Remark}
\newtheorem{Rmk}[neu]{Remark}

\newtheorem*{Ex*}{Example}

\newtheorem{Def}[neu]{Definition}
\newtheorem{Conv}[neu]{Convention}

\newcommand{\Z}{\mathbb{Z}}
\newcommand{\R}{\mathbb{R}}
\newcommand{\C}{\mathbb{C}}

\newcommand{\pf}{\longrightarrow}

\newcommand{\Mas}{\mu_{\mathrm{Maslov}}}

\newcommand{\id}{\mathrm{id}}

\newcommand{\om}{\omega}

\newcommand{\A}{\mathcal{A}}

\renewcommand{\P}{\mathcal{P}}
\newcommand{\Pe}{\mathcal{P}^{\mathrm{ess}}}

\newcommand{\D}{\mathbb{D}}

\newcommand{\M}{\mathcal{M}}
\newcommand{\Mh}{\widehat{\mathcal{M}}}
\newcommand{\Me}{\mathcal{M}^{\mathrm{ess}}}
\newcommand{\Meh}{\widehat{\mathcal{M}}^{\mathrm{ess}}}

\renewcommand{\H}{\mathrm{H}}

\newcommand{\PSS}{\mathrm{PSS}}
\newcommand{\Ham}{\mathrm{Ham}}

\newcommand{\CF}{\mathrm{CF}}
\newcommand{\CFe}{\mathrm{CF}^\mathrm{ess}}
\newcommand{\pe}{\partial^\mathrm{ess}}
\newcommand{\HF}{\mathrm{HF}}
\newcommand{\HFe}{\mathrm{HF}^\mathrm{ess}}
\newcommand{\HFl}{\mathrm{HF}^\mathrm{loc}}
\newcommand{\CM}{\mathrm{CM}}

\newcommand{\Crit}{\mathrm{Crit}}

\newcommand{\beq}{\begin{equation}}
\newcommand{\beqn}{\begin{equation}\nonumber}
\newcommand{\eeq}{\end{equation}}

\newcommand{\bea}{\begin{equation}\begin{aligned}}
\newcommand{\bean}{\begin{equation}\begin{aligned}\nonumber}
\newcommand{\eea}{\end{aligned}\end{equation}}

\numberwithin{equation}{section}

\begin{document}
\title{A note on Local Floer Homology}
\author{Peter Albers}
\address{Peter Albers\\\;Courant Institute\\\;New York University\vspace*{.3ex}}
\email{\;\htmladdnormallink{albers@cims.nyu.edu}{mailto:albers@cims.nyu.edu}\vspace*{.3ex}}
\urladdr{\;\htmladdnormallink{http://www.cims.nyu.edu/$\sim$albers}{http://www.cims.nyu.edu/~albers}}
\date{June 2006}
\subjclass[2000]{53D40, 53D12}
\begin{abstract}
In general, Lagrangian Floer homology $\HF_*(L,\phi_H(L))$ -- if well-defined -- is not isomorphic to the singular homology of the Lagrangian
submanifold $L$. For arbitrary closed Lagrangian submanifolds a local version of Floer homology $\HFl_*(L,\phi_H(L))$ is defined
in \cite{Floer_symplectic_fixed_points_and_holomorphic_spheres,Oh_Floer_cohomology_spectral_sequences_and_the_Maslov_class}
which is isomorphic to singular homology. This construction assumes that the Hamiltonian function $H$ is sufficiently $C^2$-small and the
almost complex structure involved is sufficiently standard.

In this note we develop a new construction of local Floer homology which works for any (compatible) almost complex structure and all
Hamiltonian function with Hofer norm less than the minimal (symplectic) area of a holomorphic disk or sphere. The example $S^1\subset\C$ shows that this
is sharp. If the Lagrangian submanifold is monotone, the grading of local Floer homology can be improved to a $\Z$-grading.
\end{abstract}
\maketitle
\tableofcontents
%
%
\section{Introduction}
In general, Lagrangian Floer homology
$\HF_*(L,\phi_H(L))$ -- if well-defined -- is not isomorphic to the singular homology of the Lagrangian submanifold $L$.
In \cite{Floer_symplectic_fixed_points_and_holomorphic_spheres} and \cite[section 3]{Oh_Floer_cohomology_spectral_sequences_and_the_Maslov_class}
a local version of Floer homology $\HFl_*(L,\phi_H(L))$ has been developed which is isomorphic to singular homology. This construction
assumes that the Hamiltonian function $H$ is sufficiently $C^2$-small and the almost complex structure involved is sufficiently standard. Under these
assumptions an isolating neighborhood of $L$ exists and only Floer trajectories staying inside the isolating neighborhood are considered.
In other words local Floer homology is a $C^2$-small perturbation of Morse homology.

In this note we develop a new construction of local Floer homology which works for any (compatible) almost complex structure and all Hamiltonian function
with Hofer norm less than the minimal area of a holomorphic disk or sphere. Under these much weaker assumptions an isolating
neighborhood does not exist, in general. Instead we find a concrete and \emph{geometric criterion} to single out an appropriate set of Floer
trajectories. Moreover, this criterion enables us to give direct compactness proofs for the moduli spaces involved in the construction.

The assumptions on the Hamiltonian function for this constructive approach to Floer homology are optimal. Furthermore, they place
local Floer homology is the realm of Hofer's geometry on the group of Hamiltonian diffeomorphisms.
\begin{Thm}\label{thm:thm_intro}
Let $(M,\om)$ be a closed symplectic manifold and $L\subset M$ a closed Lagrangian submanifold. We consider the Hamiltonian function
$H:S^1\times M\pf\R$ and the compatible almost structure $J$ on $(M,\om)$. The minimal area of a holomorphic sphere in $M$ or
a holomorphic disk with boundary on $L$ is denoted by $A_L$.

If the Hofer norm $||H||$ of $H$ satisfies
\beq
||H||=\int_0^1\big[\max_M H(t,\cdot)-\min_M H(t,\cdot)\big]\,dt<A_L\,,
\eeq
then there exists a distinguished subset $\Pe_L(H)$ of the set $\P_L(H)$ of Hamiltonian cords \emph{(}see equation \eqref{def:Hamiltonian_cords}\emph{)}
and for $x,y\in\Pe_L(H)$ a distinguished subset $\Me_L(x,y)$ of the moduli space of connecting Floer trajectories with the following properties.
\begin{enumerate}
\item All moduli spaces $\Me_L(x,y)$, regardless of their dimension, are compact up to breaking along elements in $\Pe_L(H)$.\\[-2ex]
\item The homology $\HFe_*(L,\phi_H(L))$ of the complex generated by $\Pe_L(H)$ with differential defined by counting zero-dimensional components
of $\Me_L(x,y)$ is canonically isomorphic to the singular homology of the Lagrangian submanifold $L$.\\[-2ex]
\item If $H_0$ and $H_1$ are two Hamiltonian functions such that $||H_0||+||H_1||<A_L$ then the well-known construction of
continuation homomorphisms carries over and provides an isomorphism $\HFe_*(L,\phi_{H_0}(L))\cong\HFe_*(L,\phi_{H_1}(L))$.
\end{enumerate}
If the Lagrangian submanifold is \emph{monotone} then $\HFe_*(L,\phi_H(L))$ is $\Z$-graded, in general it carries only a grading modulo
the minimal Maslov number of $L$.
\end{Thm}
\begin{Rmk}$ $
\begin{itemize}
\item We denote by $d_\mathrm{H}(\cdot,\cdot)$ Hofer's metric on the group $\Ham(M,\om)$ of Hamiltonian diffeomorphisms. The statement of the theorem
should be read as follows: If $\phi\in\Ham(M,\om)$ satisfies $d_\mathrm{H}(\phi,\id_M)<A_L$ then local Floer homology for $L$ and $\phi(L)$
is well defined.
\item The assumption $||H||<A_L$ is sharp, as the example $S^1\subset\C$ shows. That is, there exist Hamiltonian functions $H$ with Hofer norm
exceeding $A_{S^1}$ such that the time-1-map $\phi_H$ displace $S^1$ from itself: $S^1\cap\phi_H(S^1)=\emptyset$, in particular, $\P_L(H)=\emptyset$.
\item It is apparent from the definition of Hofer's norm that there are Hamiltonian functions $H$ with arbitrarily large $C^2$-norm
satisfying $||H||<A_L$.
\item Theorem \ref{thm:thm_intro} immediately recovers Chekanov's theorem
\cite{Chekanov_Lagrangian_intersections_symplectic_energy_and_areas_of_holomorphic_curves}
asserting that (i) the displacement energy $e(L)$ of the Lagrangian submanifold $L$ is at least as big as the
minimal area of a holomorphic disk or sphere: $e(L)\geq A_L$ and (ii), in case $||H||<A_L$ the intersection $L\cap\phi_H(L)\not=\emptyset$ contains
at least $\sum_i b_i(L)$ elements, where $b_i(L)$ are the Betti numbers of $L$.
\item The set $\Pe_L(H)$ is defined explicitly (see definition \ref{def:essential_cords}).
\end{itemize}
\end{Rmk}
Let us very briefly sketch the construction of local Floer homology according to
\cite{Floer_symplectic_fixed_points_and_holomorphic_spheres} and \cite[section 3]{Oh_Floer_cohomology_spectral_sequences_and_the_Maslov_class}.
In these articles it is proved that for sufficiently small Hamiltonian perturbations of the  Lagrangian submanifold $L$ there exists a so-called
isolating neighborhood $U$, which gives rise to a clear distinction of the set of
perturbed holomorphic strips. We recall that counting zero-dimensional families of such strips defines the boundary operator $\partial_F$ in the Floer
complex. The distinction of strips is based on the fact that either strips stay inside a compact subset of $U$ or they leave the closure $\overline{U}$.
This leads to a definition of a new boundary operator by counting only those perturbed strips which lie inside the neighborhood $U$ of $L$. In
\cite{Floer_symplectic_fixed_points_and_holomorphic_spheres,Oh_Floer_cohomology_spectral_sequences_and_the_Maslov_class} it is
proved that the new boundary operator is well-defined and the homology of the new complex equals the singular homology of $L$.

The construction of local Floer homology of a Lagrangian submanifold $L$ relies on the existence of an isolated neighborhood. The existence of such a
neighborhood is proved in the aforementioned articles (arguing by contradiction) for \emph{sufficient $C^2$-small} Hamiltonian function and for
compatible almost structures which are \emph{sufficient $C^1$-close} to the Levi-Civita almost complex structure defined in a Weinstein
neighborhood of $L$.

In section \ref{section:local_FH} we construct local Floer homology under the sole assumption that the Hofer norm of the
Hamiltonian function $H$ is less than the minimal energy $A_L$ of a holomorphic disk or sphere, that is, $||H||<A_L$ (and without any
further requirements for the compatible almost complex structure).

We specify a subset $\Pe_L(H)\subset\P_L(H)$ of the set of Hamiltonian cords. Moreover, for $x,y\in\Pe_L(H)$ we define a subset $\Me_L(x,y)$
of the moduli spaces $\M_L(x,y;J,H)$ of perturbed holomorphic strips (cf.~section \ref{section:Recollection_of_Lagrangian_Floer_homology}).
The moduli spaces $\Me_L(x,y)$ are
compact (up to breaking) given that $x,y\in\Pe_L(H)$.  Let us point out that \emph{a priori} it is not clear (to us) whether the set $\Pe_L(H)$
actually is non-empty. The above theorem proves \emph{a posteriori} that $\#\Pe_L(H)\geq\sum_i b_i(L)$.

The set $\Pe_L(H)$ of (homologically) essential cords then is used to define a new chain complex $(\CFe_*(L,\phi_H(L)),\pe_F)$ in exactly
the same way as in the usual construction of the Floer complex, namely $\CFe_*(L,\phi_H(L)):=\Pe_L(H)\otimes\Z/2$ and the differential
$\pe_F$ is defined by counting zero-dimensional components of $\Me_L(x,y)$. This results into local Floer homology $\HFe_*(L,\phi_H(L))$.
The techniques from Piunikhin, Salamon and Schwarz in \cite{Piunikhin_Salamon_Schwarz_Symplectic_Floer_Donaldson_theory_and_quantum_cohomology}
(suitably adapted to the Lagrangian setting (cf.~\cite{Albers_A_Lagrangian_PSS_morphism})) are then used to proved that $\HFe_*(L,\phi_H(L))$ is
isomorphic to $\H_*(L)$.\\[1ex]
\textbf{Acknowledgements.}\quad
The author was financially supported by the German Research Foundation (DFG) through Priority Programm 1154
"Global Differential Geometry", grant AL 904/1-1, and by NSF Grant DMS-0102298.
\section{Preliminaries}
%
\subsection{Lagrangian Floer homology}  \label{section:Recollection_of_Lagrangian_Floer_homology}
%
$ $\\[1ex]
We briefly recall the construction of the Floer complex $(\CF_i(L,\phi_H(L)),\partial_F)$ for a closed, monotone Lagrangian submanifold $L$ in
a closed, symplectic manifold $(M,\om)$.
\begin{Def}\label{def:minimal_Maslov_and_Chern_nb}
A Lagrangian submanifold $L$ of the symplectic manifolds $(M,\om)$ is called \emph{monotone}, if there exists a constant $\lambda>0$, such that
$\om|_{\pi_2(M,L)}=\lambda\cdot\Mas|_{\pi_2(M,L)}$, where $\Mas:\pi_2(M,L)\pf\Z$ is the \emph{Maslov index}.

We define the \textit{minimal Maslov number} $N_L$ of $L$ as the positive generator of the image of the Maslov index
$\Mas(\pi_2(M,L))\subset\Z$. We set $N_L=+\infty$ in case $\Mas$ vanishes. The
\emph{minimal Chern number} $N_M$ of $M$ is defined analogously.
Furthermore, we denote by $A_L$ the \emph{minimal area} of a non-constant holomorphic disk with boundary on $L$ or of a
non-constant holomorphic sphere in $M$, where area refers to the integral of the symplectic from $\om$ over the disk resp.~sphere.
\end{Def}
For a Hamiltonian function $H:S^1\times M\pf\R$ the Floer complex $\big(\CF_i(L,\phi_H(L)),\partial_F\big)$ is generated over $\Z/2$ by the set
\beq\label{def:Hamiltonian_cords}
\P_L(H):=\Big\{x\in C^{\infty}([0,1],M)\;\big|\;\dot{x}(t)=X_H\big(t,x(t)\big),\;x(0),x(1)\in L,\;[x]=0\in\pi_1(M,L)\Big\}
\eeq
i.e.~$\CF(L,\phi_H(L))=\P_L(H)\otimes\Z/2$. Let us explain some notions. First, $X_H$ is the (time dependent) Hamiltonian vector field generated
by the Hamiltonian function $H:S^1\times M\pf\R$ and is defined by $\om(X_H(t,\cdot),\,\cdot\,)=-dH(t,\cdot)$. The time-1-map
of the flow $\phi_H^t$ of $X_H$ is denoted by $\phi_H\equiv\phi_H^1$.

(A certain subset of) the intersection points $L\cap\phi_H(L)$ is often taken to generate a chain complex.
There is a one-to-one correspondence between $\P_L(H)$ and this subset by applying the flow $\phi_H^t$
to an intersection point. Furthermore, in either approach the Hamiltonian
function $H$ is required to be non-degenerate meaning that $L\pitchfork\phi_H(L)$.

The Maslov index defines a grading on $\P_L(H)$, which is only defined modulo the minimal Maslov number $N_L$ and
up to an overall shift. Let us briefly recall the construction of the grading.

Given two elements $x,y\in\P_L(H)$ we choose a map $u:[0,1]^2\rightarrow M$ s.t.~$u(0,t)=x(t)$, $u(1,t)=y(t)$ and $u(\tau,0),u(\tau,1)\in L$. According
to \cite{Viterbo_Maslov_index_for_Lagrangian_Floer_homology,Floer_A_relative_Morse_index_for_the_symplectic_action} a Maslov index is assigned to
the map $u$ as follows. Since the symplectic vector bundle $u^*TM$ is trivial, a loop of Lagrangian subspaces in $\R^{2n}$ is obtained
by following the Lagrangian subspaces of $TL$ along the two $u(\tau,0/1)$-sides of the strip and transporting them by the Hamiltonian flow along
the $u(0/1,t)$-sides (and flipping them by 90 degrees at the corners). To this loop of Lagrangian subspaces in $\R^{2n}$ the classical
Maslov index is assigned.

This gives rise to a \emph{relative} Maslov index for $x$ and $y$ which certainly depends on the choice of $u$.
Indeed, let $v:[0,1]^2\rightarrow M$ be another choice connecting $x$ and $y$ and let
$h:\D^2_+\rightarrow M$ be a half-disk realizing a homotopy of the cord $x$ to a constant path. We can form the disk $d:=h\#u\#(-v)\#(-h)$
with boundary on $L$, where $\#$ denotes concatenation and $-v$ is the map $(\tau,t)\mapsto v(-\tau,t)$. If the relative Maslov index of $x$ and $y$
is computed with help of either $u$ or $v$, the difference is given by $\Mas([d])$. We note that the Maslov index of $[d]$ does not
depend on the choice of the half-disk $h$.
Thus, we can assign a number $\mu(x,y)\in\Z/N_L$ to each pair $x,y\in\P_L(H)$. By construction, this number satisfies $\mu(x,z)=\mu(x,y)+\mu(y,z)$
for all $x,y,z\in\P_L(H)$. We artificially set $\mu(x_0):=0$ for a fixed $x_0\in\P(H)$ and define the degree $\mu(y):=\mu(y,x_0)\in\Z/N_L$
for all other $y\in\P_L(H)$. Assigning index zero to another element in $\P_L(H)$ leads to a shift of the degree. Therefore, by this procedure we define
a mod $N_L$ grading on $\P_L(H)$ up to an overall shift.

In what follows we fix the shifting ambiguity.
\begin{Def}\label{def:moduli_space_of_cut_off_holo_strips}
For $x\in\P_L(H)$ we set
\beq
\M(H;x):=\left\{d_x:\R\times[0,1]\pf M \left|\;\;
    \begin{aligned}
    &\partial_sd_x+J(t,d_x)\big(\partial_td_x-\beta(s)X_H(t,d_x)\big)=0\\
    &d_x(s,0),d_x(s,1)\in L\;\;\forall s\in\R\\
    &d_x(+\infty)=x,\;\;E(d_x)<+\infty
    \end{aligned}
    \right.
\right\}
\eeq
where $\beta:\R\rightarrow[0,1]$ is a smooth cut-off function satisfying $\beta(s)=0$ for $s\leq0$ and $\beta(s)=1$ for $s\geq1$ and
$J(t,\cdot)$, $t\in[0,1]$, is a family of compatible almost complex structures on $M$. Analogously, we set
\beq
\M(x;H):=\left\{e_x:\R\times[0,1]\pf M \left|\;\;
    \begin{aligned}
    &\partial_se_x+J(t,e_x)\big(\partial_te_x-\beta(-s)X_H(t,e_x)\big)=0\\
    &e_x(s,0),e_x(s,1)\in L\;\;\forall s\in\R\\
    &e_x(-\infty)=x,\;\;E(e_x)<+\infty
    \end{aligned}
    \right.
\right\}
\eeq
\end{Def}
By standard arguments in Floer theory it is easy to show that for generic choices of the Hamiltonian function and the almost complex structure
the moduli spaces $\M(H;x)$ and $\M(x;H)$ are smooth manifolds. Moreover, since a solution $d_x$ respectively $e_x$ has finite energy, by removal of
singularity there exists an continuous extension $d_x(-\infty)$ and $e_x(+\infty),$ respectively.

To fix the shifting ambiguity of the grading $\mu$ we require that the dimension of the moduli spaces $\M(H;x)$
is given by $\mu(x)$ mod $N_L$. Equivalently, we could demand that the space $\M(x;H)$ has dimension $n-\mu(x)$ mod $N_L$.
This convention is consistent by a gluing argument and additivity of the Fredholm index.

The Floer differential $\partial_F$ is defined by counting perturbed holomorphic strips (a.k.a.~semi-tubes or Floer strips).
For $x,y\in\P_L(H)$ we define the moduli spaces
\beq
\M_L(x,y;J,H):=\left\{\;u:\R\times[0,1]\pf M\left|\;\;
    \begin{aligned}
    &\partial_su+J(t,u)\big(\partial_tu-X_H(t,u)\big)=0\\
    &u(s,0),u(s,1)\in L\;\;\forall s\in\R\\
    &u(-\infty)=x,\;\;u(+\infty)=y
    \end{aligned}\;
\right.
\right\}
\eeq
If we would use the intersection point $L\pitchfork\phi_H(L)$
to generate the Floer complex then the differential would be defined by counting unperturbed holomorphic strips having one boundary
component on $L$ and the other on $\phi_H(L)$. Again the flow of the Hamiltonian vector field provides a one-to-one correspondence between perturbed
and unperturbed strips.

\begin{Thm}[Floer]
For a generic family $J$, the moduli spaces $\M_L(x,y;J,H)$ are smooth manifolds of dimension $\dim \M_L(x,y;J,H)\equiv\mu(y)-\mu(x)$ mod $N_L$, carrying
a free $\R$-action if $x\not=y$.
\end{Thm}

We note that the dimension of the moduli spaces is given by the Maslov index modulo the minimal Maslov number $N_L$. In other words, if we fix the
asymptotic data to be $x,y\in\P_L(H)$, the moduli space $\M_L(x,y;J,H)$ consists (in general) out of several connected components each of which has
dimension $\equiv\mu(y)-\mu(x)$ mod $N_L$.

\begin{Conv}
We set $\M_L(x,y;J,H)_{[d]}$ to be the union of the $d$-dimensional components.
\end{Conv}

\begin{Thm}[Floer, Oh]\label{thm:compactification_for_Floer_strips}
If the minimal Maslov number satisfies $N_L\geq2$ then for all $x,z\in\P_L(H)$ the moduli space
\beq
\Mh_L(x,z;J,H)_{[d-1]}:=\M_L(x,z;J,H)_{[d]}/\R
\eeq
is compact if $d=1$ and compact up to simple breaking if $d=2$, i.e.~it admits a compactification
(denoted by the same symbol) such that the boundary decomposes as follows
\beq
\partial\Mh_L(x,z;J,H)_{[1]}=\bigcup_{\substack{y\in\P_L(H)}}\Mh_L(x,y;J,H)_{[0]}\times\Mh_L(y,z;J,H)_{[0]}\,.
\eeq
\end{Thm}
\noindent
The boundary operator $\partial_F$ in the Floer complex is defined on generators $y\in\P_L(H)$ by
\beq
\partial_F(y):=\sum_{x\in\P_L(H)}\#_2\Mh_L(x,y;J,H)_{[0]}\cdot x
\eeq
and is extended linearly to $\CF_*(L,\phi_H(L))$.
Here, $\#_2\Mh_L(x,y;J,H)_{[0]}$ denotes the (mod 2) number of elements in $\Mh_L(x,y;J,H)_{[0]}$.
The two theorems above justify this definition of $\partial_F$, namely the sum is finite and $\partial_F\circ\partial_F=0$. The
\textit{Lagrangian Floer homology} groups are $\HF_*(L,\phi_H(L)):=\H_*(\CF(L,\phi_H(L),\partial_F)$.

It is an important feature of Floer homology that it is independent of the chosen family of
almost complex structure and invariant under Hamiltonian perturbations.
In particular, there exists an canonical isomorphism $\HF_*(L,\phi_H(L))\cong\HF_*(L,\phi_K(L))$ for any two Hamiltonian functions $H,K$.\\[0.5ex]
Floer theory is a (relative) Morse theory for the \emph{action functional} $\A_H$ defined on the space of paths in $M$ which start and end
on $L$ and are homotopic (relative $L$) to a constant path in $L$. By definition the action functional is
\beq\label{eqn:action_functional}
\A_H(x,d_x):=\int_{\D^2_+}d_x^*\om-\int_0^1H(t,x(t))dt
\eeq
where $d_x:\D^2_+\rightarrow M$ realizes a homotopy from a constant path to the path $x$. The value of the action functional depends only on the
relative homotopy class of $d_x$. Its critical points are exactly $\P_L(H)$.\\[2ex]
We close with a brief remark about the coefficient ring $\Z/2$. In certain cases it is possible to choose $\Z$ as coefficient ring, e.g.~if
the Lagrangian submanifold is relative spin, cf.~\cite{FOOO}. We will not pursue this direction is the present version of this article. The
same applies to non-compact symplectic manifolds which are convex at infinity or geometrically bounded.
\subsection{Some energy estimates}  \label{section:energy estimates}
%
$ $\\[1ex]
In this section we recall some standard energy estimates for elements in various moduli spaces. The derivations are simple calculations
which are carried out in \cite[appendix A]{Albers_A_Lagrangian_PSS_morphism} using the present notation.
We recall that the energy $E(u)$ of a map $u:\R\times[0,1]\pf M$ is defined as
\beq
E(u)=\!\int_{-\infty}^\infty\int_0^1|\partial_su|^2\,dt\,ds\;.
\eeq
\begin{Lemma}\label{lemma:energy_equality_for_Floer_strips}
For a Floer strip $u\in\M_L(x,y;J,H)$ the equality
\beq
E(u)=\A_H(y,d_x\#u)-\A_H(x,d_x)
\eeq
holds, where $d_x\#u$ denotes the concatenation of the half-disk $d_x$ with the Floer strip $u$.
See equation \eqref{eqn:action_functional} for the definition of the action functional.
\end{Lemma}
\noindent
For convenience we set
\beq\label{eqn:abbreviations_for_sup_and_inf}
\sup_MH:=\int_0^1\sup_MH(t,\cdot)\,dt\qquad\text{and}\qquad\inf_MH:=\int_0^1\inf_MH(t,\cdot)\,dt\;.
\eeq
In particular, in this notation the Hofer norm of a Hamiltonian function $H:S^1\times M\pf\R$ reads $\displaystyle||H||=\sup_MH-\inf_MH$.
For elements in the moduli space $\M(H;x)$ and $\M(y;H)$ from definition \ref{def:moduli_space_of_cut_off_holo_strips} we obtain the
following inequalities.
\begin{Lemma}\label{lemma:energy_inequality_for_M(H;x)}
For a solution $d_x\in\M(H;x)$ the following inequality holds:
\beq
0\leq E(d_x)\leq\;\A_H(x,d_x)+\sup_MH
\eeq
where $d_x$ serves as a homotopy from the constant path to the cord $x$.\\[1ex]
For an element $e_y\in\M(y;H)$
\beq
0\leq E(u)\leq-\A_H(y,-e_y)-\inf_MH\,,
\eeq
where $-e_y$ denotes the map $(s,t)\mapsto e_y(-s,t)$. We will actually use the following (slightly weaker) inequalities later
\beq\label{eqn:action_estimate_for_cut_off_strips}
\A_H(y,-e_y)\leq-\inf_MH\qquad\text{and}\qquad-\A_H(x,d_x)\leq\sup_MH\;.
\eeq
\end{Lemma}
\section{Local Floer homology}\label{section:local_FH}
%
\subsection{The construction}\label{section:local_FH_the_construction}
%
$ $\\[1ex]
Though some of the following makes sense for general Hamiltonian function $H$ from now on we will make the
\begin{center}
\fbox{{\textsc{Standing Assumption:}}\hspace{2ex} $\displaystyle||H||<A_L$}\\[2ex]
\end{center}
where $A_L$ is the minimal energy of an holomorphic disk or sphere and $||H||$ the Hofer norm.
\begin{Def}\label{def:essential_cords}
For a non-degenerate Hamiltonian function $H:S^1\times M\pf\R$ we set
\beq
\Pe_L(H):=\Big\{x\in\P_L(H)\mid\exists\,d_x\in\M(H;x),\;\exists\,e_x\in\M(x;H)\text{ s.t. }\om(d_x\#e_x)=0\Big\}\,,
\eeq
where the moduli spaces $\M(H;x)$ and $\M(x;H)$ are defined in \ref{def:moduli_space_of_cut_off_holo_strips}. The concatenation
of the two half disks $d_x$ and $e_x$ is denoted by $d_x\#e_x$. The integral
of the symplectic form $\om$ over this disk is denoted by $\om(d_x\#e_x)$. We call $\Pe_L(H)$ the set of (homologically) \emph{essential cords}.
\end{Def}

\begin{Rmk}
\emph{A priori} it is unclear whether $\Pe_L(H)$ is non-empty.
\end{Rmk}

\begin{Lemma}\label{lemma:homotopy_uniqueness_ess_cords}
Suppose that for $x\in\Pe_L(H)$ there exists $d_x,d_x'\in\M(H;x)$ and $e_x,e_x'\in\M(x;H)$ s.t.~$\om(d_x\#e_x)=0$ and $\om(d_x'\#e_x')=0$, then
this hold also for mixed terms: $\om(d_x\#e_x')=0$ and $\om(d_x'\#e_x)=0$ and thus
\beq
\A_H(x,d_x)=\A_H(x,d_x')=\A_H(x,-e_x)\qquad\forall d_x,d_x'\in\M(H;x)\text{ and }e_x\in\M(x;H)\,.
\eeq
\end{Lemma}

\begin{proof}
Recall that
\beq
\A_H(x,d_x)=\int_{\D^2_+}d_x^*\om-\int_0^1H(t,x(t))dt\;.
\eeq
The assumption $\om(d_x\#e_x)=0$ implies
\beq
\A_H(x,d_x)=\A_H(x,-e_x)
\eeq
where we recall that $-e_x$ is the map $(s,t)\mapsto e_x(-s,t)$. From the inequalities \eqref{eqn:action_estimate_for_cut_off_strips} we conclude
\beq\label{eqn:energy_estimate_-sup_leq_A_H(d_x)=A_H(-e_x)_leq_-inf}
-\sup_MH\leq\A_H(x,d_x)=\A_H(x,-e_x)\leq-\inf_MH\;.
\eeq
The same holds for $d_x'$ and $e_x'$
\beq
-\sup_MH\leq\A_H(x,d_x')=\A_H(x,-e_x')\leq-\inf_MH\;.
\eeq
Taking the difference of the two inequalities we obtain
\beq
\om(d_x'\#e_x)=\A_H(x,d_x')-\A_H(x,-e_x)\leq -\inf_MH+\sup_MH=||H||<A_L\;.
\eeq
On the other hand we note that by a simple gluing argument in the homotopy class $[d_x\#e_x]\in\pi_2(M,L)$ lies a holomorphic disk. Indeed
by gluing the two solutions $d_x'$, $e_x$ of Floer's equation along $x$ and then removing the Hamiltonian term a holomorphic disk is obtained, see
\cite[section 4.2.2]{Albers_A_Lagrangian_PSS_morphism} for details. In particular, $\om(d_x'\#e_x)<A_L$ implies $\om(d_x'\#e_x)=0$
by definition of $A_L$. In the same way $\om(d_x\#e_x')=0$ is proved and this immediately implies the other two equation.
\end{proof}

\begin{Rmk}
The lemma implies that we could have defined the space $\Pe_L(H)$ of essential cords by requiring that \textit{for all}
$d_x\in\M(H;x),\,e_x\in\M(x;H)$ we have $\om(d_x\#e_x)=0$.
\end{Rmk}

\begin{Def}\label{def:essential_Floer_strips}
For $x,y\in\Pe_L(H)$ we define
\beq
\Me_L(x,y):=\{u\in\M_L(x,y;J,H)\mid \om(d_x\#u\#e_y)=0\}\,,
\eeq
the set of \textit{essential Floer strips}. $\Me_L(x,y)_{[d]}$ denotes the union of the $d$-dimensional components.
\end{Def}

\begin{Rmk}$ $
\begin{enumerate}
\item By lemma \ref{lemma:homotopy_uniqueness_ess_cords} the property $\om(d_x\#u\#e_y)=0$ of an essential Floer strip does not depend on the choice of
$d_x$ or $e_y$.\\[-2ex]
\item In case $x=y\in\Pe_L(H)$ we have $\Me_L(x,x)=\M_L(x,x;J,H)=\{x\}$.\\[-2ex]
\item In case $x\not=y$ the moduli space $\Me_L(x,y)$ carries a free $\R$-action on . The quotient is denoted by
$\Meh_L(x,y)$.\\[-2ex]
\item The moduli spaces $\Me_L(x,y)$ are only defined for essential cords: $x,y\in\Pe_L(H)$. We will not mention this always but implicitly assume
that the cords are essential when we write down $\Me_L(x,y)$.
\end{enumerate}
\end{Rmk}

\begin{Def}\label{def:essential_differential}
With help of essential Floer strips we can define the differential
\beq
\pe_F(y):=\sum_{x\in\Pe_L(H)}\#_2\Meh_L(x,y)_{[0]}\cdot x\,,
\eeq
where the sum is taken over essential cords.
\end{Def}


\begin{Prop}\label{prop:boundary_operator_in_essential_Floer_homology_does_good}
For $x,z\in\Pe_L(H)$ the moduli space $\Meh_L(x,z)_{[d]}$ of essential Floer strips 
%
is compact if $d=0$ and compact up to simple breaking along essential cords if $d=1$, i.e.~it admits a compactification
(denoted by the same symbol) such that the boundary decomposes as follows
\beq
\partial\Meh_L(x,z)_{[1]}=\bigcup_{\substack{y\in\Pe_L(H)}}\Meh_L(x,y)_{[0]}\times\Meh_L(y,z)_{[0]}\,.
\eeq
We note, that the union is taken over essential cords.
\end{Prop}

\begin{proof}
We start with a simple observation which actually was the starting point of this approach to local Floer homology.
For $u\in\Me_L(x,y)$ the following energy estimate holds (and is proved below)
\beq\label{eqn:energy_estimate_essential_strips}
E(u)<A_L\;.
\eeq
Indeed, recall from lemma \ref{lemma:energy_equality_for_Floer_strips}
\beq
E(u)=\A_H(y,d_x\#u)-\A_H(x,d_x)\;.
\eeq
Since $u$ is essential, $\om(d_x\#u\#e_y)=0$ holds, i.e.~$\om(d_x\#u)=\om(-e_y)$. This implies
$\A_H(y,d_x\#u)=\A_H(y,-e_y)$. Now we can apply inequality \eqref{eqn:action_estimate_for_cut_off_strips} from
lemma \ref{lemma:energy_inequality_for_M(H;x)} to conclude
\beq
E(u)=\A_H(y,-e_y)-\A_H(x,d_x)\leq-\inf_MH+\sup_MH=||H||\;.
\eeq
According to our standing assumption we obtain equation \eqref{eqn:energy_estimate_essential_strips}. This allows
us to prove that the moduli spaces $\Meh_L(x,z)$ are compact up to breaking.
Let us assume that a sequence $(u_n)\subset\Me_L(x,y)$ develops a bubble, for instance $u_n$ converges in the Gromov-Hausdorff topology to
$(u_\infty,D)$, where $D$ is a holomorphic disk, then (cf.~\cite[proposition 4.6.1]{McDuff_Salamon_J_holomorphic_curves_and_symplectic_topology}
and \cite[proposition 3.3]{Salamon_lectures_on_floer_homology})
\beq
E(u_\infty)+E(D)\leq E(u_n)\leq||H||<A_L\;.
\eeq
This immediately implies that $E(D)<A_L$ and thus the holomorphic disk $D$ is constant. This obviously generalizes to multiple bubbling of holomorphic
disks and spheres. In particular, all moduli spaces $\Me_L(x,y)_{[d]}$ are compact up to breaking for all dimensions $d$. \\[1ex]
To finish the proof of the proposition we need to show that breaking occurs only along essential cords and broken Floer strips are essential.
Let us assume that we have a sequence $(u_n)\subset\Me_L(x,z)$ which converges to a broken solution
$(v_1,v_2)\in\M_L(x,y;J,H)\times\M_L(y,z;J,H)$.\\[1.5ex]
We are required to prove that $y\in\Pe_L(H)$ and $v_1\in\Me_L(x,y)$ and $v_2\in\Me_L(y,z)$.\\[1.5ex]
Pick $d_x\in\M(H;x)$ and $e_z\in\M(z;H)$. By gluing $d_x$ and $v_1$ we find an element $d_y\in\M(H;y)$ and, in turn, by gluing $v_2$ and $e_z$
we find an element $e_y\in\M(y;H)$. The gluing is the standard gluing of
two Floer strips. Since the homotopy class is preserved in the limit $u_n\rightharpoonup(v_1,v_2)$, i.e.~$\om(u_n)=\om(v_1)+\om(v_2)$, we derive
\beq
\om(d_y\#e_y)=\om(d_x\#v_1\#v_2\#e_z)=\om(d_x\#u_n\#e_z)=0\;,
\eeq
because $u_n$ is essential.
In particular, $y$ is an essential cord $y\in\Pe_L(H)$. Moreover, the Floer strips $v_1$ and $v_2$ are essential, since
\beqn
\om(d_x\#v_1\#e_y)=\om(d_x\#v_1\#v_2\#e_z)=0 \quad\text{and}\quad \om(d_y\#v_2\#e_z)=\om(d_x\#v_1\#v_2\#e_z)=0\;.
\eeq
This concludes the proof of the proposition.
\end{proof}
\begin{Rmk}
Proposition \ref{prop:boundary_operator_in_essential_Floer_homology_does_good}
shows that $\pe_F$ is well-defined and a differential: $\pe_F\circ\pe_F=0$. We note that we do not use any monotonicity assumption for the Lagrangian
submanifold $L$.

Moreover, we proved more, namely all moduli space $\Me_L(x,y)_{[d]}$ are compact up to breaking regardless of their dimension $d$ and they can be
compactified by essential Floer strips.
\end{Rmk}
\begin{Def}\label{def:essential_Floer_homology}
We set $\displaystyle \CFe_*(L,\phi_H(L)):=\Pe_L(H)\otimes\Z/2$ and
\beq
\HFe_*(L,\phi_H(L)):=\H_*(\CFe(L,\phi_H(L)),\partial_F)\;.
\eeq
We call $\HFe_*(L,\phi_H(L))$ \emph{local Floer homology}.
\end{Def}

So far the grading of $\CFe_*(L,\phi_H(L))$ is, as described in section \ref{section:Recollection_of_Lagrangian_Floer_homology},
a $\Z/N_L$-grading plus an overall shifting ambiguity.

\begin{Lemma}
If the Lagrangian submanifold $L$ is \textit{monotone}, i.e.~$\om|_{\pi_2(M,L)}=\lambda\cdot\Mas$ for some $\lambda>0$,
then the $\Z/N_L$-grading of $\CFe_*(L,\phi_H(L))$ can be improved to a $\Z$-grading (still with shifting ambiguity).
Moreover, the differential $\pe_F$ preserves the $\Z$-grading (and not only the $\Z/N_L$-grading).
\end{Lemma}

\begin{proof}
The grading on the Floer complex is defined by assigning a Maslov index to pairs of cords $x,y\in\P_L(H)$ (see section
\ref{section:Recollection_of_Lagrangian_Floer_homology}). This involves the choice of a map $u:[0,1]^2\pf M$ with the
properties $u(0,t)=x(t)$, $u(1,t)=y(t)$ and $u(\tau,0),u(\tau,1)\in L$. Different homotopy classes of such maps change the Maslov index
by multiples of the minimal Maslov number $N_L$. Thus, a $\Z/N_L$-grading is obtained.

For essential cords $x,y\in\Pe_L(H)$ of a Hamiltonian function $H$ satisfying $||H||<A_L$ there is a preferred choice
of a (homotopy class of a) map $u:[0,1]^2\pf M$, namely such that $\om(d_x\#u\#e_y)=0$. We recall that the maps $d_x\in\M(H;x)$ and $e_x\in\M(x;H)$
exists by definition of essential cords (see definition \ref{def:essential_cords}).

Because of the monotonicity of $L$ we claim that for all choices of such a map $u$ the relative Maslov index for the pair $x,y\in\Pe_L(H)$ give
rise to the same value. Indeed, let us assume that we choose two maps $u,v$ satisfy $\om(d_x\#u\#e_y)=0$ and $\om(d_x\#v\#e_y)=0$.
The difference of the relative Maslov index computed with $u$ or $v$ is given by the Maslov index of the disk
$D:=d_x\#u\#(-v)\#(-d_x)$. Under the assumption that $L$ is monotone we compute
\bea
\lambda\cdot\Mas(D)&=\om(d_x\#u\#(-v)\#(-d_x))\\
        &=\om(d_x\#u\#e_y)+\om((-e_y)\#(-v)\#(-d_x))\\
        &=\om(d_x\#u\#e_y)-\om(d_x\#v\#e_y)\\
        &=0
\eea
In particular, we can compute the relative Maslov index of $x$ and $y$ with help of $u$ or $v$ equally well.
The differential $\pe_F$ is defined by using essential Floer strips $u$ i.e.~$\om(d_x\#u\#e_y)=0$, thus,
$\pe_F$ preserves the $\Z$-grading.
\end{proof}

\begin{Conv}\label{conv:grading}
As proved above, if $L$ is monotone we obtain a $\Z$-grading for local Floer homology, but in general only a $\Z/N_L$-grading.
All subsequent statements have to be read accordingly.
\end{Conv}

\begin{Prop}
For two Hamiltonian functions $H_0$ and $H_1$ satisfying
\beq
||H_0||+||H_1||< A_L
\eeq
the (obvious modification of the) continuation homomorphisms are well-defined and provide isomorphism between the local Floer homologies of
$H_0$ and $H_1$.
\end{Prop}

\begin{proof}
We consider the homotopy $H_s:=\beta(s)H_1+(1-\beta(s))H_0$ where $\beta(s)$ is a smooth cut-off function satisfying $\beta(s)=0$ for $s\leq0$
and $\beta(s)=1$ for $s\geq1$. To define the continuation homomorphisms in Floer homology
the homotopy parameter $s$ is coupled to the $\R$-parameter in the Floer equation, cf.~e.g.~\cite[section 3.4]{Salamon_lectures_on_floer_homology}.
That is, counting solutions $u$ of $\partial_su+J(t,u)\big(\partial_tu-X_{H_s}(t,u)\big)=0$ with $u(-\infty)=x\in\Pe_L(H_0)$ and
$u(+\infty)=y\in\Pe_L(H_1)$ defines the continuation homomorphism. Furthermore, we require that $\om(d_x\#u\#e_y)=0$. Then the energy estimate
from lemma \ref{lemma:energy_equality_for_Floer_strips} changes
into
\beq
E(u)\leq \A_{H_1}(y,d_x\#u)-\A_{H_0}(x,d_x)+\int_0^1\sup_M\big[H_1(t,\cdot)-H_0(t,\cdot)\big]dt\;.
\eeq
Therefore, as in the proof of proposition \ref{prop:boundary_operator_in_essential_Floer_homology_does_good} we use $\om(d_x\#u\#e_y)=0$
and the inequalities
\eqref{eqn:action_estimate_for_cut_off_strips} to conclude
\beq
E(u)\leq -\inf_MH_1+\sup_MH_0+\sup_MH_1-\inf_MH_0=||H_0||+||H_1||
\eeq
(using the notation from \eqref{eqn:abbreviations_for_sup_and_inf}). The compactness arguments as in the proof of proposition
\ref{prop:boundary_operator_in_essential_Floer_homology_does_good} carry over unchanged. Thus the appropriate moduli spaces are compact up to breaking
along essential cords and counting defines a map $\HFe_*(H_1)\pf\HFe_*(H_0)$. The inverse
is constructed by interchanging the roles of $H_0$ and $H_1$. We leave the details to the reader.
\end{proof}

In the construction of the articles
\cite{Floer_symplectic_fixed_points_and_holomorphic_spheres,Oh_Floer_cohomology_spectral_sequences_and_the_Maslov_class} the above proposition is
proved as well but again under the assumption that the homotopy is $C^2$-small and the almost complex structure is
sufficiently close to the Levi-Civita almost complex structure.

Moreover, in the mentioned articles the proposition is crucial for proving that the local Floer homology is isomorphic to the singular homology of the
Lagrangian submanifold $L$. Namely, choosing a $C^2$-small Morse function on $L$ and the Levi-Civita almost complex
structure, the local Floer complex reduces to the Morse complex of the Morse function $f$.

The construction of $\HFe(L,\phi_H(L))$ is designed in such a way that the techniques from Piunikhin, Salamon and Schwarz in
\cite{Piunikhin_Salamon_Schwarz_Symplectic_Floer_Donaldson_theory_and_quantum_cohomology} can be applied to Lagrangian Floer homology.
%
\subsection{The isomorphism}\label{section:local_FH_the_isomorphism}
%
$ $\\[1ex]
In this section we prove that $\HFe(L,\phi_H(L))$ is canonically isomorphic to the singular homology of the Lagrangian submanifold $L$.
This is achieved by applying the ideas from \cite{Piunikhin_Salamon_Schwarz_Symplectic_Floer_Donaldson_theory_and_quantum_cohomology}. In
\cite{Piunikhin_Salamon_Schwarz_Symplectic_Floer_Donaldson_theory_and_quantum_cohomology} an
isomorphism $\PSS:\H^{n-k}(M)\stackrel{\cong}{\pf}\HF_k(H)$ between the \textit{Hamiltonian Floer homology} of the Hamiltonian
function $H$ and the singular homology of the manifold $M$ for a very general class of symplectic manifolds $(M,\om)$ is constructed.

The analogous result for \textit{Lagrangian Floer homology} $\HF_*(L,\phi_H(L))$ is certainly false due to the existence of displaceable Lagrangian
submanifolds. In general, Lagrangian Floer homology is only well-defined for monotone Lagrangian submanifolds with minimal Maslov number $N_L\geq2$
(cf.~\cite{Floer_Morse_theory_for_Lagrangian_intersections,Oh_Floer_cohomology_of_Lagrangian_intersections_and_pseudo-holomorphic_disks_I}).
The question to what extend the techniques from \cite{Piunikhin_Salamon_Schwarz_Symplectic_Floer_Donaldson_theory_and_quantum_cohomology} can be
carried over to Lagrangian Floer homology is addressed in \cite{Albers_A_Lagrangian_PSS_morphism}.\\[.5ex]
\indent We recall the standing assumption $||H||<A_L$.

\subsubsection{Lagrangian Piunikhin-Salamon-Schwarz morphisms}\label{section:Lagrangian_PSS_morphisms}
%
\begin{Thm}[\cite{Albers_A_Lagrangian_PSS_morphism}, theorem 1.1]\label{thm:existence_of_Lagrangian_PSS}
We consider a $2n$-dimensional closed symplectic manifold $(M,\om)$ and a closed monotone Lagrangian
submanifold $L\subset M$ of minimal Maslov number $N_L\geq2$. Then there exist homomorphisms
\begin{gather}
\varphi_k:\HF_k(L,\phi_H(L))\pf\H^{n-k}(L;\Z/2)\quad\text{for }k\leq N_L-2\;,\\[1ex]
\rho_k:\H^{n-k}(L;\Z/2)\pf\HF_k(L,\phi_H(L))\quad\text{for }k\geq n-N_L+2\;,
\end{gather}
where $H:S^1\times M\pf\R$ is a Hamiltonian function and $\phi_H$ its time-1-map.

For $n-N_L+2\leq k\leq N_L-2$ the maps are inverse to each other
\beq
\varphi_k\circ\rho_k=\id_{\H^{n-k}(L;\Z/2)}\quad\text{and}\quad
        \rho_k\circ\varphi_k=\id_{\HF_k(L,\phi_H(L))}\;.
\eeq
\end{Thm}
\begin{Rmk}

The homomorphisms in theorem \ref{thm:existence_of_Lagrangian_PSS} are constructed using the ideas introduced by Piunikhin, Salamon and Schwarz in
\cite{Piunikhin_Salamon_Schwarz_Symplectic_Floer_Donaldson_theory_and_quantum_cohomology}. We call them \textsc{Lagrangian PSS morphisms}.
The restrictions on the degrees are sharp in general as examples show (see \cite[remark 2.6]{Albers_A_Lagrangian_PSS_morphism}).
\end{Rmk}

The Lagrangian PSS morphisms are defined by counting zero-dimensional components of certain moduli spaces. In fact, the moduli
spaces defining $\varphi$ are the intersection of the space $\M(H;x)$ with the unstable manifolds of some critical point of a Morse function on $L$. The
moduli space corresponding to $\rho$ is defined by  the intersection of $\M(x;H)$ with some stable manifold (see equation
\eqref{eqn:def_of_moduli_spaces_for_varphi_and_rho} for details). The degree restrictions in theorem
\ref{thm:existence_of_Lagrangian_PSS} for the Lagrangian PSS morphisms are due to bubbling-off of holomorphic disks, i.e.~non-compactness of
the spaces $\M(H;x)$ and $\M(x;H)$. In the present situation bubbling can be ruled out.

\begin{Prop}\label{prop:compactness_of_M(x;H)_and_M(H;x)_for_essential_cords}
The moduli spaces $\M(H;x)$ and $\M(x;H)$ (cf.~definition \ref{def:moduli_space_of_cut_off_holo_strips})
are compact up to breaking for all $x\in\Pe_L(H)$. Moreover, they can be compactified by adding essential Floer strips.
\end{Prop}

\begin{proof}
Since $x$ is essential both moduli spaces $\M(H;x)$ and $\M(x;H)$ are non-empty we recall inequality
\eqref{eqn:energy_estimate_-sup_leq_A_H(d_x)=A_H(-e_x)_leq_-inf} from the proof of lemma \ref{lemma:homotopy_uniqueness_ess_cords}
\beq
-\sup_MH\leq\A_H(x,d_x)=\A_H(x,-e_x)\leq-\inf_MH\;.
\eeq
The observation is that if both moduli spaces are non-empty each one provides compactness for the other. Indeed, for
$d_x\in\M(x;H)$ we combined the energy estimate in lemma \ref{lemma:energy_inequality_for_M(H;x)} with the inequality from above and obtain
\beq
E(d_x)\leq\A_H(x,d_x)+\sup_MH\leq-\inf_MH+\sup_MH=||H||\;.
\eeq
Analogously, for $e_x\in\M(H;x)$ we obtain
\beq
E(e_x)\leq-\A_H(x,-e_x)-\inf_MH\leq\sup_MH-\inf_MH=||H||\;.
\eeq
By assumption, the Hamiltonian function $H$ satisfies $||H||<A_L$ and we can argue as in the proof of proposition
\ref{prop:boundary_operator_in_essential_Floer_homology_does_good} to exclude bubbling-off and to conclude that we need only to add
essential Floer strips in the compactification.
\end{proof}

\begin{Rmk}
Again we point out that the above compactness result for the moduli spaces $\M(H;x)$ and $\M(x;H)$ holds regardless of their dimension as
long as $x$ is essential.
\end{Rmk}

%
\subsubsection{The proof of $\;\HFe_*(L,\phi_H(L))\cong\H_*(L)$}\label{section:Lagrangian_PSS_morphisms_local}
%
$ $\\[1ex]
The next theorem is the adaptation of theorem \ref{thm:existence_of_Lagrangian_PSS} to local Floer homology.
\begin{Thm}\label{thm:existence_of_Lagrangian_PSS_LOCAL_case}
Let $L$ be a closed (not necessarily monotone) Lagrangian submanifold $L$ in a closed symplectic manifold $(M,\om)$. If the
Hamiltonian function $H$ satisfies
\beq
||H||<A_L\,,
\eeq
then there exist isomorphisms
\begin{gather}
\varphi_*:\HFe_*(L,\phi_H(L))\pf\H^{n-*}(L;\Z/2)\;,\\[1ex]
\rho_*:\H^{n-*}(L;\Z/2)\pf\HFe_*(L,\phi_H(L))\;,
\end{gather}
moreover, \quad$\varphi_*=\rho_*^{-1}$.
\end{Thm}

\begin{Rmk}
The homomorphism $\varphi_*$ and $\rho_*$ are the restrictions of the Lagrangian PSS morphisms (appearing in theorem
\ref{thm:existence_of_Lagrangian_PSS}) to the set of
essential Hamiltonian cords. We denote them by the same symbol.
\end{Rmk}

\begin{proof}
The proof is basically the same as the in \cite{Albers_A_Lagrangian_PSS_morphism} though greatly simplified due to the standing assumption $||H||<A_L$.
Let us briefly recall the definition of the maps $\varphi_*$ and $\rho_*$.
For full details see \cite[section 4.1]{Albers_A_Lagrangian_PSS_morphism}. As usual the maps are defined by a counting process.
\bea\label{eqn:def_of_moduli_spaces_for_varphi_and_rho}
\M^{\varphi}(q,x)&:=\M(H;x)\pitchfork_{\,\mathrm{ev}}\!W^u(q;f)\;,\\[1ex]
\M^{\rho}(x,p)&:=\M(x;H)\pitchfork_{\,\mathrm{ev}}\!W^s(q;f)\;.
\eea
Let us explain the notation. Recall (cf.~definition \ref{def:moduli_space_of_cut_off_holo_strips})
that elements $u\in\M(H;x)$ admit a continuous extension $u(-\infty)$. The evaluation map $\mathrm{ev}:\M(H;x)\pf L$
assigns to each element $u$ this value: $\mathrm{ev}(u)=u(-\infty)$. The moduli space $\M^{\varphi}(q,x)$ consists out of those maps $u$ for
which $u(-\infty)$ lie in the unstable manifold $W^u(q;f)$ of the critical point $q$ of a Morse function $f:L\pf\R$.
Furthermore, we want the evaluation map $\mathrm{ev}$ to be
transverse to the unstable manifolds. Analogously, the second definition has to be read, where the maps in $\M(x;H)$ are evaluated at $+\infty$.

We denote elements in the moduli space $\M^{\varphi}(q,x)$ as pairs $(\gamma,u)$ where $u\in\M(H;x)$ and $\gamma:(-\infty,0]\pf L$ is a solution
of $\dot{\gamma}(t)=-\nabla f(\gamma(t))$ and $\gamma(-\infty)=q$.

It is straight forward to prove that for generic choices these are smooth manifolds. But in general bubbling-off can occur for sequences
in these moduli spaces. The compactness properties are governed by those of the moduli spaces $\M(x;H)$ and $\M(H;x)$. In the general case, this leads
to the restrictions appearing in theorem \ref{thm:existence_of_Lagrangian_PSS}.

The standing assumption $||H||<A_L$ provides compactness (up to breaking) of the moduli spaces $\M(x;H)$ and $\M(H;x)$ (cf.~proposition
\ref{prop:compactness_of_M(x;H)_and_M(H;x)_for_essential_cords}). In particular, the same is true for $\M^{\varphi}(q,x)$ and $\M^{\rho}(x,q)$. Thus,
the maps $\varphi_*$ and $\rho_*$ are well-defined for all degrees.
\bea
\varphi_k:\CFe_k(L,\phi_H(L))&\pf\CM^{n-k}(f)\\
                    x\quad&\mapsto \sum_{q\in\Crit(f)}\#_2\M^{\varphi}(q,x)_{[0]}\cdot q
\eea

\bea
\rho_k:\CM^{n-k}(f)&\pf\CFe_k(L,\phi_H(L))\\
                    q\quad&\mapsto \sum_{x\in\Pe_L(H)}\#_2\M^{\rho}(x,q)_{[0]}\cdot q
\eea
Since the moduli spaces $\M^{\varphi}(q,x)$ and $\M^{\rho}(x,q)$ are compact up to breaking along essential cords these maps descent to homology.\\[1ex]
More involved is to prove that they are inverse to each other. We need to consider the compositions $\varphi_*\circ\rho_*$ and
$\rho_*\circ\varphi_*$. The idea in both cases is to find a suitable cobordism relating the composition to the identity map.

The easier case is $\rho_*\circ\varphi_*$  since it again relies only on the compactness (up to breaking) of the moduli spaces
$\M(H;x)$ and $\M(x;H)$. The geometric idea is as follows.
\begin{enumerate}
\item The composition $\rho_k\circ\varphi_k$ is a map from Floer homology to Floer homology. The coefficient
of $\rho_k\circ\varphi_k(y)$ in front of $x\in\P_L(H)$, it is given by counting all zero dimensional configurations
$(u_+,\gamma_+;\gamma_-,u_-)$, where $(u_+,\gamma_+)\in\M^{\rho}(x,q)$ and $(\gamma_-,u_-)\in\M^{\varphi}(q,y)$ and $q\in\Crit(f)$ is arbitrary.

\item We glue the two gradient flow half-trajectories $\gamma_+$ and $\gamma_-$ at the critical point $q$ and obtain $(u_+,\Gamma,u_-)$, where
$\Gamma$ is a \emph{finite length} gradient flow trajectory, say parameterized by $[0,R]$, such that $u_+(+\infty)=\Gamma(0)$ and
$\Gamma(R)=u_-(-\infty)$.

\item We shrink the length $R$ of the gradient flow trajectory to zero. In the limit $R=0$ we obtain a pair $(u_+,u_-)$ of two maps
$u_+,u_-:\R\times[0,1]\pf M$ which satisfy Floer's equation on one half and are holomorphic on the $\mp$-half of the strip. Furthermore,
they satisfy $u_+(+\infty)=u_-(-\infty)$ and $u_+(-\infty)=x$ and $u_-(+\infty)=y$.

\item Since $(u_+,u_-)$ both are holomorphic around the point $u_+(+\infty)=u_-(-\infty)$ we can glue them and obtain a map
$u:\R\times[0,1]\pf M$ which satisfies Floer's equation with Hamiltonian term given by $H$ up to a compact perturbation around $s=0$.

\item We remove the compact perturbation and obtain an honest Floer strip connecting $x$ to $y$. Since we count zero dimensional configurations
(and are not dividing out the $\R$-action),
this is only non-zero if $y=x$, in which case there is exactly one such strip, namely the constant.

In other words, up to a cobordism (i.e.~in homology), the coefficient $\rho_k\circ\varphi_k(y)$
in front of $x$ equals zero or one depending on whether $x=y$.
\end{enumerate}

For a detailed definition of the moduli spaces involved as well as for a series of pictures illustrating the idea
we refer the reader to \cite[section 4.2.1]{Albers_A_Lagrangian_PSS_morphism}. From the
above description it is apparent that only the compactness (up to breaking) of the moduli spaces $\M(H;x)$ and $\M(x;H)$ is an issue. Since
we assume that $||H||<A_L$ this poses no problem here. Indeed, proposition \ref{prop:compactness_of_M(x;H)_and_M(H;x)_for_essential_cords}
guarantees compactness for the moduli spaces $\M(H;x)$ and $\M(x;H)$ in all dimensions.\\[1ex]
The more delicate composition to handle is $\varphi_*\circ\rho_*$. Again we sketch the idea (see
\cite[section 4.2.2]{Albers_A_Lagrangian_PSS_morphism} for some pictures).
\begin{enumerate}
\item The coefficient of $\varphi_k\circ\rho_k(p)$ in front of $q\in\Crit(f)$ is given by counting zero dimensional
configurations $(\gamma_-,u_-;u_+,\gamma_+)$ such that $(\gamma_-,u_-)\in\M^{\varphi}(q,x)$ and $(u_+,\gamma_+)\in\M^{\rho}(x,p)$
for some $x\in\P_L(H)$.

\item We glue $u_-$ and $u_+$ at the cord $x\in\P_L(H)$ and obtain a single strip $U:\R\times[0,1]\pf M$ which is a solution of Floer's equation.
The important fact to note is, that the Hamiltonian term in the Floer equation is zero outside a compact subset of $\R\times[0,1]$
and that $U$ satisfies $\gamma_-(0)=U(-\infty)$ and $U(+\infty)=\gamma_+(0)$.

We note that the Morse indices of $q$ and $p$ are equal. The set of triples $(\gamma_-,U,\gamma_+)$ as described above
is obtained by intersecting the space of maps $U$ with the unstable manifold of $q$ and the stable manifold of $p$. In particular, the space formed by
the maps $U$ has to be of dimension $n=\dim L$. This implies that the Maslov index on the (relative) homotopy class $[U]\in\pi_2(M,L)$ of $U$ equals
zero: $\Mas([U])=0$. Furthermore, since $x$ is essential we conclude $\om(U)=0$.

\item The compact perturbation by the Hamiltonian term can be removed and we end up with triples $(\gamma_-,U,\gamma_+)$, where $U$ is a holomorphic
map $U:\R\times[0,1]\pf M$ (of finite energy) satisfying $\gamma_-(0)=U(-\infty)$ and $U(+\infty)=\gamma_+(0)$. Thus, $U$ is a \emph{holomorphic
disk} with boundary on the Lagrangian submanifold $L$.

\item The integral $\om(U)$ vanishes and therefore, $U$ has to be constant and $(\gamma_-,\gamma_+)$ form
an gradient flow line from $q$ to $p$. Again we are interested in zero dimensional configuration and we are not dividing by the $\R$-action.
By the same arguments as before we obtain the identity map $\id_{\H^{n-k}(L;\Z/2)}$.
\end{enumerate}
The problem is that we need to consider a new kind of moduli space which is not of the types considered so far. Let us be more precise.

For $q,p\in\Crit(f)$ we define the moduli space $\M^{\varphi\circ\rho}(q,p)$ to be the set of quadruples $(R,\gamma_-,U,\gamma_+)$, where
\beq
R\geq0,\quad\gamma_-:(-\infty,0]\pf L,\quad U:\R\times[0,1]\pf M,\quad\gamma_+:[0,+\infty)\pf L
\eeq
satisfying
\begin{gather}
\partial_sU+J(t,U)\big(\partial_tU-\tilde{\alpha}_R(s)\cdot X_H(t,U)\big)=0\,,\\[1.5ex]
U(s,0),\,U(s,1)\in L,\quad E(U)<+\infty\,,\\[1.5ex]
\dot{\gamma}_\pm(t)+\nabla^g f\circ\gamma_\pm(t)=0\,,\\[1.5ex]
\gamma_-(-\infty)=q,\quad\gamma_-(0)=U(-\infty),\quad U(+\infty)=\gamma_+(0),\quad\gamma_+(+\infty)=p\;.
\end{gather}
Finally, we demand that the relative homotopy class $[U]\in\pi_2(M,L)$ satiesfies
\beq\label{eqn:homotopy_class=0_for_rho_o_phi}
\Mas([U])=0\quad\text{and}\quad\om([U])=0\;.
\eeq
The map $\tilde{\alpha}_R$ is a cut-off function such that for $R\geq1$ we have $\tilde{\alpha}_R(s)=1$ for $|s|\leq R$ and $\tilde{\alpha}_R(s)=0$
for $|s|\geq R+1$. Furthermore, we require for its slope that $-1\leq\tilde{\alpha}_R'(s)\leq0$ for $s\geq0$ and $0\leq\tilde{\alpha}_R'(s)\leq1$ for
$s\leq0$. For $R\leq1$ we set
$\tilde{\alpha}_R(s)=R\tilde{\alpha}_1(s)$. In particular, for $R=0$ the cut-off function vanishes identically: $\tilde{\alpha}_0\equiv0$.

That this space is a smooth manifold for generic choices is a again achieved by standard arguments.
Of course, compactness problems are only caused by sequences $(U_N)$ of maps $U_N:\R\times[0,1]\pf M$ from above.

The following energy estimate is easily derived (cf.~\cite[Lemma A.3]{Albers_A_Lagrangian_PSS_morphism})
\beq
0\leq E(U_n)\leq\om([U_n])+||H||\;.
\eeq
Since we require $\om([U_n])=0$ and by our standing assumption we conclude $E(U_n)<A_L$. In particular, following the arguments
in the proof of proposition \ref{prop:boundary_operator_in_essential_Floer_homology_does_good}, the moduli spaces
$\dim\M^{\varphi\circ\rho}(q,p)$ are compact up to breaking again for all dimensions. Without the assumption $||H||<A_L$
theorem \ref{thm:existence_of_Lagrangian_PSS} this is not true, in general.

Since the moduli spaces $\M^{\varphi\circ\rho}(q,p)$ are compact up to breaking counting zero-dimensional components defines a map
$\Theta^{\varphi\circ\rho}_k:\CM^{n-k}(L;\Z/2)\pf\CM^{n-k-1}(L;\Z/2)$. From the compactification of the one-dimensional components of
$\M^{\varphi\circ\rho}(q,p)$ it is apparent that $\Theta^{\varphi\circ\rho}$ is a chain homotopy between $\varphi\circ\rho$ and the identity.

All further details can be found in \cite[section 4.2.2]{Albers_A_Lagrangian_PSS_morphism}, in particular in the proof of theorem 4.25.
\end{proof}

As an immediate corollary of theorem \ref{thm:existence_of_Lagrangian_PSS_LOCAL_case} we obtain Chekanov's result.
\begin{Cor}[\cite{Chekanov_Lagrangian_intersections_symplectic_energy_and_areas_of_holomorphic_curves}]\label{cor:chekanov}
Let $L$ be a closed Lagrangian submanifold in a closed symplectic manifold $(M,\om)$. Denote by $A_L$ the minimal energy of a holomorphic disk
with boundary on $L$ or a holomorphic sphere. If the Hamiltonian function $H:S^1\times M\pf\R$ is non-degenerate and has Hofer norm less than $A_L$,
\beq
||H||<A_L\,,
\eeq
then $\#\P_L(H)\geq\#\Pe_L(H)\geq\sum_ib_i(L)$. In particular, there exists at least $\sum_ib_i(L)$ Hamiltonian cords with action bounded as follows
\beq\label{eqn:action_estimate_for_essential_strips}
-\sup_MH\leq\A_H(x,d_x)\leq-\inf_MH\;.
\eeq
\end{Cor}

\begin{Rmk}
Chekanov proves this result for all geometrically bounded symplectic manifolds. The methods used in this article carry over to this case.
Moreover, it seems that the action bounds \eqref{eqn:action_estimate_for_essential_strips}
cannot be derived from Chekanov's approach since he uses some abstract homological algebra.
\end{Rmk}
\noindent\hrulefill

%
%
\bibliographystyle{amsalpha}
\bibliography{../../../Bibtex/bibtex_paper_list}
\end{document}